\documentclass{ACmod}

\DeclareFontFamily{U}{rsf}{}
\DeclareFontShape{U}{rsf}{m}{n}{
  <5> <6> rsfs5 <7> <8> <9> rsfs7 <10-> rsfs10}{}
\DeclareMathAlphabet{\mathscr}{U}{rsf}{m}{n}

\DeclareMathAlphabet{\mathgth}{U}{euf}{m}{n}

\DeclareFontFamily{U}{cyr}{}
\DeclareFontShape{U}{cyr}{m}{n}{
  <5> wncyr5 <6> wncyr6 <7> wncyr7 <8> wncyr8 <9> wncyr9 <10-> wncyr10}{}
\DeclareMathAlphabet{\mathcyr}{U}{cyr}{m}{n}

\input cyracc.def

\DeclareSymbolFont{bbold}{U}{bbold}{m}{n}
\DeclareSymbolFontAlphabet{\mathbbold}{bbold}

\makeatletter
\def\operator@font{\sf}
\makeatother

\newcommand{\cA}{{\mathscr A}}

\newcommand{\cE}{{\mathscr E}}

\newcommand{\cU}{{\mathscr U}}

\newcommand{\cV}{{\mathscr V}}

\newcommand{\MC}{{\sf MC}}

\newcommand{\MF}{{\sf MF}}

\newcommand{\Jac}{{\mathsf{Jac}}}
\newcommand{\hres}{{\mathsf{hres}}}

\newcommand{\bone}{{\mathbf 1}}

\DeclareMathOperator{\Tr}{Tr}

\DeclareMathOperator{\Der}{Der}

\DeclareMathOperator{\End}{End}

\DeclareMathOperator{\SL}{SL}

\DeclareMathOperator{\id}{id}

\DeclareMathOperator{\KS}{KS}

\newcommand{\ra}{\rightarrow}

\newcommand{\C}{\mathbb{C}}

\newcommand{\Z}{\mathbb{Z}}

\renewcommand{\phi}{\varphi}

\usepackage{tikz}
\usepackage{amscd}
\usepackage[tocflat]{tocstyle}
\usetikzlibrary{positioning}
\usetikzlibrary{shapes.misc}
\numberwithin{equation}{section}

\tikzset{cross/.style={cross out, draw=black, minimum size=2*(#1-\pgflinewidth), inner sep=0pt, outer sep=0pt},
cross/.default={1pt}}

\begin{document}

\title{Categorical Saito theory, I: A comparison result}

\author[Junwu Tu]{Junwu Tu\footnote{Partially supported by the NSF grant DMS-1801806. } }

\address{Junwu Tu, Institute of Mathematical Science, ShanghaiTech University, Shanghai 201210, China. Email: tujw@shanghaitech.edu.cn}

\begin{abstract} 
  {\sc Abstract:} In this paper, we present an explicit cyclic minimal $A_\infty$ model for the category of matrix factorizations $\MF(W)$ of an isolated hypersurface singularity. The key observation is to use Kontsevich's deformation quantization technique. Pushing this idea further, we use the Tsygan formality map to obtain a comparison theorem that the categorical Variation of Semi-infinite Hodge Structure of $\MF(W)$ is isomorphic to Saito's original geometric construction in primitive form theory.  An immediate corollary of this comparison result is that the analogue of C\u ald\u araru's conjecture holds for the category $\MF(W)$. 
\end{abstract}

 \maketitle

\setcounter{tocdepth}{1}

\section{Introduction}

The notion of Variation of Semi-infinite Hodge Structures (VSHS) was first introduced by Barannikov~\cite{Bar}, inspired by earlier works of Saito~\cite{Sai2} and~\cite{Sai1}. Roughly speaking, a (formal) VSHS over a complete regular local ring $R=\mathbb{K}[[t_1,\cdots,t_\mu]]$ consists of a $R[[u]]$-module $\cE$, locally free of finite rank, a flat meromorphic connection $\nabla$ in both $t$- and $u$-directions, and a $R$-linear bilinear form $\langle -,-\rangle: \cE\otimes_R \cE \rightarrow R[[u]]$.
The triple $(\cE,\nabla,\langle-,-\rangle)$ must satisfy certain compatibility conditions, detailed in Section~\ref{sec:vshs}. Let $S:=\mathbb{K}[[x_1,\cdots,x_N]]$ be the formal power series ring, and let $W\in S$ be a function with an isolated singularity at the origin. The main result of this paper is a comparison theorem between two VSHS's associated with the function $W$: one constructed by Saito~\cite{Sai2} known as primitive form theory, and the other one from non-commutative Hodge theory of the category of matrix factorizations $\MF(W)$~\cite{KKP},~\cite{Shk},~\cite{Shk2},~\cite{GPS}. In particular, this comparison result implies that there is a natural bijection between categorical primitive forms of $\MF(W)$ and Saito's primitive forms of $W$, and that the corresponding Frobenius manifolds obtained are also isomorphic. In a sequel~\cite{Tu}, we shall give a categorical construction of genus zero $B$-model invariants of Landau-Ginzburg orbifolds, conjecturally mirror to the genus zero part of the FJRW theory. It is not yet known how to extend Saito's geometric construction to the orbifold case.

\paragraph{{\bf Saito's VSHS.}}  Denote by
\[ \Omega_S^*:= \mathbb{K}[[x_1,\cdots,x_N]][dx_1,\cdots,dx_N]\]
the space of differential forms of the commutative algebra $S$. Since $W$ has an isolated singularity, its Jacobian ring $\Jac(W):=\mathbb{K}[[x_1,\cdots,x_N]]/(\partial_1W,\cdots,\partial_NW)$ is finite dimensional.  
We may choose a basis $\phi_1,\cdots,\phi_\mu$ of $\Jac(W)$. Saito~\cite{Sai2} defined a VSHS on the twisted de Rham complex of the formal mini-versal deformation $\mathscr{W}:=W+t_1\phi_1+\cdots+t_\mu\phi_\mu$ of $W$. More precisely, this VSHS is given by
\[\cV^{{\sf Saito}}:=\Big( H^*\big( \Omega_S^*[[t_1,\cdots,t_\mu,u]], d\mathscr{W}+ud_{DR} \big), \nabla, K^\mathscr{W}\Big),\] 
with $\nabla_{t_j}:=\frac{d}{d t_j} - \frac{\phi_j}{u}$, $\nabla_{u}:=\frac{d}{du} -\frac{N\cdot \id}{2u}+ \frac{\mathscr{W}}{u^2}$~\footnote{This differs from Saito's original definition by the term $-\frac{N\cdot\id}{2u}$.}, and $K^\mathscr{W}$ is the higher residue pairing extending the residue pairing.

\paragraph{{\bf Categorical construction of VSHS's.}} It has long been suspected, for example in Kontsevich-Manin~\cite{KonMan}, Saito-Takahashi~\cite{SaiTak}, Shklyarov~\cite{Shk3}~\cite{Shk4} that Saito's construction may be put into a categorical setting, where one replaces geometric deformation of singularities by deformations of differential graded (or more generally $A_\infty$) categories. In this paper, we shall achieve this through Kontsevich's deformation quantization technique. Indeed, we consider the category $\MF(W)$ of matrix factorizations of $W$. An important structure result of Dyckerhoff~\cite{Dyc} asserts that the category $\MF(W)$ is compactly generated by a single object denoted by $\mathbb{K}^{{\sf stab}}$, the stabilization of residue field. This reduces the study of deformations of the category $\MF(W)$ to deformations of the differential graded algebra $\End(\mathbb{K}^{{\sf stab}})$. In fact, we shall pass to its minimal model $A_\infty$ algebra
\[ A^W:= \mbox{ the minimal model of }\; \End(\mathbb{K}^{{\sf stab}}).\]
Next, we may use the deformed potential function $\mathscr{W}=W+t_1\phi_1+\cdots+t_\mu\phi_\mu$ to obtain a formal deformation $A^{\mathscr{W}}$ of the $A_\infty$ algebra $A^W$. The non-commutative VSHS associated to $\MF(W)$ is then given by
\[ \cV^{{\MF(W)}}:= \Big( HC^-_*(A^{\mathscr{W}}), \nabla, \langle-,-\rangle_{{\sf hres}}\Big),\]
with the connection in the $t$-direction given by Getzler's connection~\cite{Get}, and in the $u$-direction the connection operator defined in~\cite{KonSoi}\cite{KKP}\cite{Shk}\cite{CLT}. The pairing $\langle-,-\rangle_{{\sf hres}}$ is the categorical higher residue pairing defined in~\cite{Shk2}\cite{She}. Note that the definition of $\cV^{{\MF(W)}}$ is independent of the choice of the compact generator $\mathbb{K}^{{\sf stab}}$, as all the structures involved are Morita-invariant~\cite{She}.

\paragraph{{\bf Deformation quantizations.}}  The main result of this paper states that there is an isomorphism of VSHS's between $\cV^{{\sf Saito}}$ and $\cV^{{\MF(W)}}$. To obtain such a comparison result, we need to use Kontsevich's deformation quantization technique. 

To begin with, recall Dyckerhoff~\cite{Dyc} showed that the minimal model $A_\infty$ algebra is, as a vector space, given by
\[ A^W \cong \mathbb{K}[\epsilon_1,\cdots,\epsilon_N],\]
the super-commutative algebra in $N$ odd variables. To obtain the $A_\infty$ algebra structure on the vector space, the standard tool is using homological perturbation technique and Konsevich-Soibelman's tree formulas. Implementing this idea, Dyckerhoff obtained some partial formulas of the higher products. The key observation of this paper is to instead use Kontsevich's deformation quantization to obtain such a minimal model. Indeed, consider the super-commutative algebra 
\[A:=\mathbb{K}[\epsilon_1,\cdots,\epsilon_N],\]
and its Koszul dual algebra $S=\mathbb{K}[[x_1,\cdots,x_N]]$ with $x_j$ linear dual variables of $\epsilon_j$.  Associated to the algebra $A$, we form the space of poly-vector fields:
\begin{equation}~\label{eq:poly-vector}
 T_{{\sf poly}} (A)=\mathbb{K}[\epsilon_1,\cdots,\epsilon_N]\otimes \mathbb{K}[[ \partial_{\epsilon_1},\cdots,\partial{\epsilon_N}]],
 \end{equation}
with $\epsilon_j$'s degree $1$ and $\partial_{\epsilon_j}$'s degree $0$. The shifted space $T_{{\sf poly}} (A)[1]$, endowed with zero differential and the Schouten bracket forms a differential graded Lie algebras (DGLA). Now, given an element $W(x_1,\cdots,x_N)\in S$, we may consider it as a poly-vector fields with constant coefficients
\[ W(\partial_{\epsilon_1},\cdots,\partial_{\epsilon_N})\in T_{{\sf poly}}(A)[1].\]
Observe that since it has constant coefficients, it satisfies the Maurer-Cartan equation
\[ [W(\partial_{\epsilon_1},\cdots,\partial_{\epsilon_N}), W(\partial_{\epsilon_1},\cdots,\partial_{\epsilon_N})]=0.\]

At this point, we envoke Kontsevich's formality morphism
\[ \mathcal{U}: T_{{\sf poly}} (A)[1] \ra C^*(A)[1],\]
an $L_\infty$ quasi-isomorphism from the DGLA of poly-vector fields to the shifted Hochschild cochain complex which is endowed with the Hochschild differential and Gerstenhaber bracket. The push-forward $\mathcal{U}_*W$ is thus a Maurer-Cartan element of $C^*(A)[1]$ which by definition gives an $A_\infty$ structure on the vector space $A$. In Corollary~\ref{cor:minimal}, we prove that this $A_\infty$ structure is quasi-isomorphic to Dyckerhoff's dga $\End(\mathbb{K}^{{\sf stab}})$. 

\paragraph{{\bf Cyclicity of the minimal model $A^W$.}} We single out one of many nice features of Kontsevich's formula: the minimal model $A^W$ is automatically cyclic, with respect to the trace map
\[ \Tr: \mathbb{K}[ \epsilon_1,\ldots,\epsilon_N ]  \ra \mathbb{K}, \;\;\; \Tr(\epsilon_I)=\begin{cases} 1, \mbox{\;\;\; if $I=(1,2,3,\cdots,N)$}\\ 0, \mbox{\;\;\; otherwise.} \end{cases}\]
This follows from a result of Felder-Shoikhet~\cite{FelSho} (which deals with the divergence-free case, see Willwatcher-Calaque~\cite{WilCal} for the general case). We explain why having a cyclic $A_\infty$ model is important for us. Indeed, our original motivation was to perform explicit computation of Costello's categorical Gromov-Witten invariants for the category of matrix factorizations, where one needs a cyclic $A_\infty$ algebra as the input data. Previously, the only known (non-trivial) example of cyclic $A_\infty$ structure is due to Polishchuk in the case of elliptic curves. This was used in~\cite{CalTu} to perform computation of Costello's invariant in the case of $g=1$, $n=1$. Using Kontsevich's explicit formula, we obtain a cyclic $A_\infty$ structure for every potential $W\in S$! In particular, this can be applied to perform explicit computation of categorical Gromov-Witten invariants for the derived category of a Quintic, using Orlov's equivalence known as the Calabi-Yau/Landau-Ginzburg correspondence. This line of research will be pursued elsewhere.

\paragraph{{\bf Comparison via Tsygan formality map.}} Tsygan formality asserts that there exists an $L_\infty$ module quasi-isomorphism
\[ \mathcal{U}^{Sh}: CC_*(A) \rightarrow \Omega^*_A.\]
Here the superscript ``Sh" is for Shoikhet, following Willwacher~\cite{Wil}. Let $W\in T_{{\sf poly}}(A)$ be a poly-vector field with constant coefficients. Then the Tsygan formality map can be deformed by $W$. The first deformed structure map $\mathcal{U}^{Sh,W}_0$ is explicitly given by
 \[ \mathcal{U}^{Sh,W}_0(-)=\sum_{k\geq 0} \frac{1}{k!} \mathcal{U}^{Sh}_k(\underbrace{W,\cdots,W}_{\mbox{$k$-copies}};-).\]
Formally, over the ring $R=\mathbb{K}[[t_1,\cdots,t_\mu]]$, we may further deform the Tsygan formality map by $\mathscr{W}= W + t_1\varphi_1+\cdots+ t_\mu\varphi_\mu$
to obtain a chain map
 \[  \mathcal{U}^{Sh,\mathscr{W}}_0:  \big( CC_*(A^{\mathscr{W}}), b^{\mathscr{W}}\big) \rightarrow \big( \Omega^*_A\otimes \mathbb{K}[[t_1,\cdots,t_\mu]], L_{\mathscr{W}}\big).\]
Extend this map $u$-linearly to obtain a map which we still denote by $ \mathcal{U}^{Sh,\mathscr{W}}_0$ on the negative cyclic chain complex
\[
\mathcal{U}_0^{Sh,\mathscr{W}}: \big( CC_*(A^{\mathscr{W}})[[u]], b^{\mathscr{W}}+uB\big) \rightarrow \big( \Omega^*_A\otimes \mathbb{K}[[t_1,\cdots,t_\mu]][[u]], L_{\mathscr{W}}+ud_{DR}\big).\]
The fact that $\mathcal{U}_0^{Sh,\mathscr{W}}$ remains a chain map (i.e. it intertwines the Connes operator $B$ with the de Rham differential $d_{DR}$) is proved by Willwatcher~\cite[Theorem 1.3]{Wil}. The deformed Tsygan formality map may be considered as a B-model open-closed map. Indeed, it is defined using integrals over configuration spaces of points in the unit disk, where one inserts $\mathscr{W}$ in the interior markings, as illustrated below.
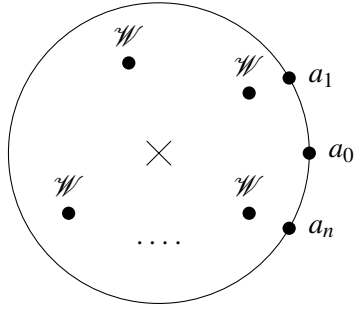
\[ \begin{tikzpicture}[baseline={(current bounding box.center)},scale=0.4]
\draw (0,0) circle (5);
\node[circle,fill=black,inner sep=0pt,minimum size=5pt,label=right:{$a_0$}] at (5,0) {};
\node[circle,fill=black,inner sep=0pt,minimum size=5pt,label=right:{$a_1$}] at (4.33,2.5) {};
\node[circle,fill=black,inner sep=0pt,minimum size=5pt,label=right:{$a_n$}] at (4.33,-2.5) {};
\draw (0,0) node[cross=5pt] {};
\node[circle,fill=black,inner sep=0pt,minimum size=5pt,label=above:{$\mathscr{W}$}] at (3,2) {};
\node[circle,fill=black,inner sep=0pt,minimum size=5pt,label=above:{$\mathscr{W}$}] at (-1,3) {};
\node[circle,fill=black,inner sep=0pt,minimum size=5pt,label=above:{$\mathscr{W}$}] at (-3,-2) {};
\node[circle,fill=black,inner sep=0pt,minimum size=5pt,label=above:{$\mathscr{W}$}] at (3,-2) {};
\node at (0,-3) {$\cdots\cdot$};
\end{tikzpicture}\]

\begin{Theorem}~\label{intro:thm}
The deformed Tsygan formality map induces an isomorphism 
\[\mathcal{U}_0^{Sh,\mathscr{W}}: HC^-_*(A^{\mathscr{W}}) \rightarrow H^*\big( \Omega^*_A[[t_1,\cdots,t_\mu,u]], L_{\mathscr{W}}+ud_{DR}\big)\]
of VSHS's, up to a constant in $\mathbb{K}^*$. (The right hand side VSHS is naturally isomorphic to Saito's $\cV^{{\sf Saito}}$ see Section~\ref{sec:tsygan}.)
\end{Theorem}

\medskip
\paragraph{{\bf C\u ad\u araru's conjecture for the category $\MF(W)$.}} In~\cite{Cal} and~\cite{CalWil}, Hochschild structures associated with the derived category of coherent sheaves on a smooth projective variety $X$ were studied in great detail. In particular, one obtains a ring structure on $HH^*(X)$, a $HH^*(X)$-module structure on $HH_*(X)$, and a symmetric non-degenerate pairing $\langle-,-\rangle_{{\sf Muk}}: HH_*(X)\otimes HH_*(X) \ra \mathbb{K}$. There is another much easier but similar structure: a ring structure on $H^*(X, \bigoplus_{j=0}^{\dim X} \Lambda^jT_X[-j])$, a module structure on $H^*(X,  \bigoplus_{j=0}^{\dim X} \Omega^j_X[j])$, and the usual Poincare pairing. It was conjectured in {\em loc. cit.} that there exist isomorphisms
\[ I^K: HH^*(X) \ra H^*(X, \bigoplus_{j=0}^{\dim X} \Lambda^jT_X[-j]), \;\;\; I_K: HH_*(X) \ra H^*(X,  \bigoplus_{j=0}^{\dim X} \Omega^j_X[j]),\]
which would intertwine all the mentioned structures. This is known as the C\u ad\u araru's conjecture, later proved by Calaque-Rossi-Van den Bergh~\cite{CRV} for the ring and module structure, and by Ramadoss~\cite{Ram} for the pairing. 

Theorem~\ref{intro:thm} above in particular implies the analogue of C\u ad\u araru's conjecture holds for the category $\MF(W)$. More precisely, we have the following
\begin{Theorem}~\label{thm:caldararu}
Assume that $W\in \mathbb{K}[[x_1,\cdots,x_N]]$ has an isolated singularity at origin. Let $d\mathcal{U}_W: {{\sf Jac}} (W) \ra HH^*(A^W)$ be the tangent map of Kontsevich's formality morphism at $W$. And let $\mathcal{U}_0^{Sh,W}: HH_*(A^W) \ra {{\sf Jac}}(W)dx_1\wedge\cdots\wedge dx_N$ be the Tsygan's formality map deformed by $W$. Then the pair $(d\mathcal{U}_W,\mathcal{U}_0^{Sh,W})$ intertwines the ring structures, the module structures, and the pairings up to a constant in $\mathbb{K}^*$. Note that on the space ${{\sf Jac}}(W)dx_1\wedge\cdots\wedge dx_N$ we use the residue pairing. 
\end{Theorem}

\begin{Proof}
The fact that the tangent map preserves the rings structure was already proved by Kontsevich in his original deformation quantization paper~\cite[Section 8]{Kon2}. The module structure is preserved follows from the definition of the Getzler connection that the $u^{-1}$-term is given by capping with the Kodaira-Spencer class, and that the Kodaira-Spencer map is an isomorphism. The map $\mathcal{U}_0^{Sh,W}$ intertwines the pairing up to constant is part of the above comparison theorem, since the higher residue pairing is an extension of the residue pairing.
\end{Proof}

\paragraph{{\bf Organization of the paper.}} In Section~\ref{sec:vshs}, we recollect the notion of VSHS's and its appearance in non-commuative Hodge theory. In Section~\ref{sec:kontsevich}, we use Kontsevich's deformation quantization formula to write down an explicit cyclic $A_\infty$ miminmal model for the category $\MF(W)$. In Section~\ref{sec:tsygan}, we prove Theorem~\ref{intro:thm} using the deformed Tsygan formality map. 

\paragraph{{\bf Conventions.}} We following Sheridan's sign convention and notations in~\cite{She}. In particular, the notation $\mu_n (n\geq 1)$ stands for higher $A_\infty$ products in the shifted sign convention. Following~\cite[Section 3]{She}, the notations $b^{1|1}$ and $B^{1|1}$ stands for certain actions of Hochschild cochains on Hochschild Chains. These operators give a calculus structure with $\iota=b^{1|1}+uB^{1|1}$ the contraction operator, and the Lie derivative action defined by
\begin{align*}
 L_\phi(a_0|a_1|\cdots |a_n):=  & \sum (-1)^{( |a_0|'+\cdots+|a_j|')\cdot |\phi|'} a_0|a_1|\cdots|a_j|\phi(a_{j+1},\cdots,a_{j+l})|\cdots|a_n\\
 &+\sum \phi(\overbrace{a_0,\cdots,a_{l-1}})|a_l|\cdots |a_n
\end{align*}

\paragraph {\bf Acknowledgments.} The author is grateful to Andrei C\u ald\u araru for many useful discussions and suggestions at various stages of the work. As already mentioned in the introduction, the paper originated from our long-term joint project on computing Costello's categorical invariants.  The author also thanks Nicholas Addington and Sasha Polishchuk for an inspiring discussion on the subject and their hospitality during a short visit at University of Oregon. Finally, special thanks to the Forschungsinstitut f\" ur Mathematik at ETH Z\" urich where the author had the great opportunity to present results of the paper  in a beautiful workshop organized by Andrei C\u ald\u araru, Rahul Pandharipande, and Nick Sheridan. 
\section{Variation of Semi-infinite Hodge structures}~\label{sec:vshs}

In this section, we recall the notion of VSHS's and its natural appearance in non-commutative Hodge theory.

\paragraph{{\bf Semi-infinite Hodge structures.}} Semi-infinite Hodge Structure (and its variations) was introduced by  Barannikov~\cite{Bar} (see also Barannikov-Kontsevich~\cite{BarKon}), based on earlier works of Saito~\cite{Sai2}~\cite{Sai1}. This structure naturally appears in the categorical setup, known as non-commutative Hodge theory~\cite{KKP},~\cite{Shk}~\cite{Shk2},~\cite{She}. First, we recall the following definition, essentially from~\cite[Section 2]{She}, the difference being that we also have a connection in the $u$-direction. We abbreviate Semi-infinite Hodge Structure by SHS.

\begin{Definition}~\label{def:shs}
A SHS of dimension $N\in \Z/2\Z$ consists of 
\begin{itemize}
\item A $\Z/2\Z$-graded $\mathbb{K}[[u]]$-module $E$ that is free and of finite rank.
\item A meromorphic connection of even degree $\nabla_{\frac{\partial}{\partial u}}: E \rightarrow u^{-2} E$, with at most a second order pole at $u=0$.
\item An even $\mathbb{K}$-linear pairing $ \langle-,-\rangle: E\otimes_\mathbb{K} E \rightarrow \mathbb{K}[[u]]$.
\end{itemize}
They are subjected to the following conditions:
\begin{itemize}
\item The pairing is $\mathbb{K}[[u]]$-sesquilinear, i.e. $\langle f(u) s_1, g(u) s_2\rangle= f(u)g(u)^{\star}\langle s_1, s_2\rangle$. Here the operator $\star: \mathbb{K}[[u]]\rightarrow \mathbb{K}[[u]]$ is given by $g(u)^\star= g(-u)$.
\item The pairing is covariantly constant with respect the $u$-connection.
\item The pairing is anti-symmetric in the sense that $\langle s_1, s_2\rangle = (-1)^{N+|s_1||s_2|}\langle s_2, s_1 \rangle^\star$.
\end{itemize}
A SHS is called polarized if the induced pairing $E/uE\otimes E/uE \rightarrow \mathbb{K}$ is non-degenerate.
\end{Definition}

\medskip
If $E$ is SHS (polarized or unpolarized), we shall refer to the decreasing filtration
\[ F^k E := u^k\cdot E, \;\; k\geq 0\]
the Hodge filtration of $E$.

\paragraph{{\bf SHS's from non-commutative geometry.}} Let $A$ be an $A_\infty$ algebra over $\mathbb{K}$. We shall always assume that $A$ is strict unital. Let $CC_*(A)$ denote its reduced Hochschild chain complex, with the Hochschild differential $b$. Its homology is denoted by $HH_*(A)$. The Hochschild chain complex is endowed with a circle action (i.e. a degree one, square zero operator commuting with the differential $b$) $B$, known as the Connes differential $B(a_0|a_1|\cdots|a_k):=\sum_j \pm 1|a_j|\cdots|a_k|a_0|\cdots|a_{j-1}$ where the sum is over cyclic permutations of the indices. From this circle action, we obtain the periodic cyclic chain complex $\big( CC_*(A)[u^{-1},u]], b+uB\big)$, the negative cyclic chain complex $\big( CC_*(A)[[u]], b+uB\big)$, and the cyclic chain complex $\big( CC_*(A)[u^{-1}] u^{-1}, b+uB\big)$. Note that the negative cyclic chain complex is a sub-complex of the periodic cyclic chain complex, while the cyclic chain complex is a quotient. Their corresponding homology groups will be denoted by $HP_*(A)$, $HC^-_*(A)$, and $HC_*(A)$ respectively.

A fundamental conjecture of noncommutative Hodge theory, as formulated by Konstevich-Soibelman~\cite{KonSoi} and Katzarkov-Kontsevich-Pantev~\cite{KKP}, is the Hodge-to-de-Rham degeneration property:

\begin{Conjecture}~\label{conj:hodge-to-de-Rham}
Let $A$ be a $\Z/2\Z$-graded, compact, smooth $A_\infty$ algebra. Then its negative cyclic homology $HC_*^-(A)$ is a locally free $\mathbb{K}[[u]]$-module of finite rank.
\end{Conjecture}

\vspace{0.5cm}

In order to obtain a SHS on $HC_*^-(A)$, we need to have a connection in the $u$-direction.
This connection was first introduced by Kontsevich-Soibelman~\cite{KonSoi}, and Shklyarov~\cite{Shk2}\cite{Shk4}. Its formula in the $A_\infty$ case was written down in~\cite{CLT}. 

We can also define a pairing on $ HC^-_*(A)$. Indeed, by~\cite[Section 5.3]{She}, there exists a chain level Mukai-pairing, i.e. a chain map of degree zero 
\[  \langle-,-\rangle_{{\sf Muk}}: CC_*(A)^{\otimes 2} \rightarrow \mathbb{K}.\]
Assuming that $A$ is compact and smooth, the non-degeneracy property of the Mukai-pairing is proved by Shklyarov~\cite{Shk}.

The chain level Mukai-pairing, extended sesquilinearly to the negative cyclic chain complex $CC_*(A)[[u]]$, induces the so-called higher residue pairing
 \begin{align*}
 \langle -,-\rangle_\hres :& CC_{*}(A)[[ u]]^{\otimes 2} \rightarrow \mathbb{K}[[u]].
 \end{align*}
It follows from the definition of the higher residue pairing that it is compatible with the cyclic differential $b+uB$. Hence it descends to the negative cyclic homology to give a $\mathbb{K}$-linear pairing
\[ \langle-,-\rangle_\hres: HC_{*}^-(A) \otimes HC_{*}^-(A) \rightarrow \mathbb{K}[[u]].\]
However, to obtain the anti-symmetry property of this pairing. We further assume that the $A_\infty$ algebra is endowed with a cyclic pairing of parity $N$. Then, the chain level Mukai pairing is symmetric after shifting the complex by $N\pmod{2}$. More precisely, we have the identity
\[ \langle \alpha, \beta \rangle_{{\sf Muk}}= (-1)^{(N+|\alpha|)(N+|\beta|)}\langle \beta, \alpha\rangle_{{\sf Muk}}.\]
The anti-symmetric property of $\langle-,-\rangle_\hres$ follows from the above symmetric property of $\langle-,-\rangle_{{\sf Muk}}$.  Putting everything  together, we arrive at the following 

\begin{Theorem}~\label{thm:shs}
Let $A$ be a $\Z/2\Z$-graded, compact and smooth cyclic $A_\infty$ algebra of parity $N\pmod{2}$. Assume Conjecture~\ref{conj:hodge-to-de-Rham} holds for $A$. Then the triple $\big( HC_*^-(A), \nabla_{\frac{\partial}{\partial u}}, \langle-,-\rangle_{\hres}\big)$ forms a polarized SHS.
\end{Theorem}

\medskip
\begin{remark}
If the $HH_*(A)$ is concentrated at either purely even or odd degree, then Conjecture~\ref{conj:hodge-to-de-Rham} trivially holds by degree reason. For the purpose of this paper, the Hochschild homology of the category $\MF(W)$ (assuming $W$ has an isolated singularity) indeed satisfies this condition.
\end{remark}

\paragraph{{\bf Variational Semi-infinite Hodge Structures (VSHS).}} Let $(R,\mathfrak{m})$ be a complete regular local ring of finite type over $\mathbb{K}$. A polarized VSHS over $R$ is given by
\begin{itemize}
\item A $\Z/2\Z$-graded free $R[[u]]$-module $\mathcal{E}$ of finite rank.
\item A {{\sl flat}}, degree zero, meromorphic connection $\nabla: \mathcal{E} \rightarrow \Omega_{R[[u]]/\mathbb{K}} \otimes_{R[[u]]} \mathcal{E}$ such that it has a simple pole along $u=0$ in the $R$-direction, and at most an order two pole at $u=0$ in the $u$-direction. In other words, we have
\[ \nabla_X: \cE \ra u^{-1} \cE, \; X\in \Der(R), \; \;\; \nabla_{\frac{\partial}{\partial u}}: \cE \ra u^{-2} \cE.\]
\item A $R$-linear, $u$-sesquilinear, and $\nabla$-constant pairing $\langle-,-\rangle: \mathcal{E}\otimes \mathcal{E} \rightarrow R[[u]]$. 
\end{itemize}
The above data satisfies the same conditions as in Definition~\ref{def:shs}. We only note that since $\mathcal{E}$ is free of finite rank, the nondegeneracy condition of
\[ \mathcal{E}/u\mathcal{E} \otimes_R \mathcal{E}/u\mathcal{E} \rightarrow R\]
is equivalent to the nondegeneracy condition at the central fiber $E=\mathcal{E}/\mathfrak{m}\mathcal{E}$, i.e. the induced pairing
\[ E/uE \otimes_\mathbb{K} E/uE \ra \mathbb{K}\]
is non-degenerate.

\paragraph{{\bf VSHS from flat families of $A_\infty$ algebras.}} Let $\cA$ be an $A_\infty$ algebras over $(R,\mathfrak{m})$, by which we mean $\cA$ is a $\Z/2\Z$-graded free $R$-module of finite rank, endowed with an $A_\infty$ structure that is $R$-linear. There are some technicalities involved in establishing a generalization of Theorem~\ref{thm:shs} to a family version. 

Indeed, we begin with a simple counter-example. Let $R=\mathbb{K}[[t]]$, and let $\cA=R\oplus R\cdot \epsilon$ a rank two, free $R$-module, with the first copy $R$ at even degree and the second copy $R$ at odd degree. Endow $\cA$ with an unital $R$-linear differential graded algebra structure:
\[ \mu_1(\epsilon)=t\cdot \bone, \;\; \mu_2(\epsilon,\epsilon)=\bone.\]
Let $CC_*(\cA)$ be its usual Hochschild chain complex. By direct computation, one can show that its homology is given by $\mathbb{K}[t^{-1}]$, with the $R$-module structure defined by
\[ t. t^{-k}:=\begin{cases} t^{1-k}, \mbox{\;\; if $k\geq 1$,}\\
0, \mbox{\;\; if $k=0$}\end{cases}\]
Thus, this is an $R$-module with $0$ generic fiber, and with an infinite rank special fiber.

In order to resolve this issue, we introduce the $\mathfrak{m}$-adic completed Hochschild chain complex. More precisely we set
\[ CC_*^{\mathfrak{m}}(\cA):= \varprojlim CC_*(\cA/\mathfrak{m}^k).\]
Since this is the only version of Hochschild chain complex we shall use in this paper,  in the following we keep the notation $CC_*(\cA)$ for its $\mathfrak{m}$-adic completed version.  

Continuing with the example above, one can show that the completed version of the Hochschild homology of $\cA$ is a free $R$-module of rank one, generated by the element
\[ \sum_{k\geq 0}   (-1)^k \frac{t^k}{2^k} \epsilon|\epsilon^k.\]
Notice that this infinite sum makes sense in the $\mathfrak{m}$-adic completed Hochschild chain complex. 

The periodic cyclic chain complex of $\cA$ is defined by the inverse limit \[\varprojlim_{k} \big( CC_*(\cA/ \mathfrak{m}^k \cA)[u^{-1},u]], b+uB\big) .\] 
Similarly, we also define the cyclic chain complex/negative cyclic chain complex as the $\mathfrak{m}$-adic completion the corresponding chain complex associated with the $A_\infty$ algebra $\cA/\mathfrak{m}^k\cA$.

The same as in the absolute case, the negative cyclic homology of $\cA$ carries the following structures:
\begin{itemize}
\item The higher residue pairing
\[ \langle-,-\rangle_{\hres}: HC^-_*(\cA) \otimes HC^-_*(\cA) \ra R[[u]].\]
\item The Getzler connection
\[ \nabla^{{{\sf Get}}}: \Der (R) \otimes HC^-_*(\cA) \ra u^{-1} HC^-_*(\cA).\]
This connection has a first order pole at $u=0$.
\item An $u$-direction connection
\[ \nabla_{\frac{\partial}{\partial u}}: HC^-_*(\cA) \ra u^{-2} HC^-_*(\cA).\]
It has a second order pole at $u=0$.
\item If $\cA$ is a formal deformation of $A$ parametrized by a regular local ring $R$, then the $u$-direction connection commutes with the Getzler connection in the $R$-direction. That is, for any vector field $\nu \in \Der (R)$, we have
\[ [\nabla_{\frac{\partial}{\partial u}}, \nabla_\nu^{{\sf Get}}]=0.\]
\end{itemize}

Thus if we set $\nabla:=\nabla^{{{\sf Get}}}+\nabla_{\frac{\partial}{\partial u}}$, it defines a flat connection both in the $R$-direction and the $u$-direction.

Finally, we recall the following definition from~\cite{CLT}.

\begin{Definition}
A formal deformation $\cA$ of $A$ over $R$, with structure maps $\mu=\prod_k \mu_k$, is called versal if the Kodaira-Spencer map
\begin{equation*}
 \KS : \Der (R) \ra HH^*(\cA), \;\; \KS(\frac{\partial}{\partial t_j}):= [\frac{\partial \mu}{\partial t_j}],
 \end{equation*}
is an isomorphism.
\end{Definition}

\begin{Lemma}
Let $A$ be a $\Z/2\Z$-graded compact, smooth, cyclic $A_\infty$ algebra, such that its Hochschild cohomology $HH^*(A)$ is concentrated at purely even degree. Let $\cA$ be its versal formal deformation over a regular local ring $(R,\mathfrak{m})$. Then for any $k\geq 1$, we have
\[ HH_*(\cA/\mathfrak{m}^k\cA)\cong  HH_*(A) \otimes_\mathbb{K} R/\mathfrak{m}^k.\]
Similarly, for the negative cyclic homology, we have 
\[ HC^-_*(\cA/\mathfrak{m}^k\cA)\cong  HC^-_*(A) \otimes_\mathbb{K} R/\mathfrak{m}^k.\]
\end{Lemma}

\begin{Proof} 
The existence of a versal formal deformation of $A$ follows from that $HH^*(A)$ is concentrated at purely even degree, so that the obstruction space of its deformation theory vanishes. We also need to use the compactness and smoothness of $A$ to ensure that $HH^*(A)$ is finite dimensional, so that the deformation space is of finite type. 

To prove the locally freeness of $HH_*(\cA/\mathfrak{m}^k\cA)$, choose a homotopy retraction $(i,p,h): HH_*(A) \cong CC_*(A)$, which exists since we are over a field $\mathbb{K}$. Extend it $R/\mathfrak{m}^k$-linearly to get a homotopy retraction which we still denote by
\[ (i,p,h): HH_*(A) \otimes_\mathbb{K} R/\mathfrak{m}^k \cong CC_*(A)\otimes_\mathbb{K} R/\mathfrak{m}^k.\]
The Hochschild complex $CC_*(\cA/\mathfrak{m}^k\cA)$ is a deformation of the right hand side complex. Indeed, its differential is given by $b^{(0)}+b^{(1)}+\cdots+b^{(k-1)}$ where $b^{(i)}$ denotes the differential induced by the $i$-th order higher multiplication of the deformation $\cA$. The perturbation $b^{(1)}+\cdots+b^{(k-1)}$ is a small perturbation since every term in it is at least first order, which implies that higher enough compositions of it will vanish over the ring $R/\mathfrak{m}^k$. Thus we obtain a perturbed homotopy retraction
\[ (i',p',h') : HH_*(A) \otimes_\mathbb{K} R/\mathfrak{m}^k \cong CC_*(\cA/\mathfrak{m}^k\cA).\]
In general, the differential on the left hand side should also be perturbed accordingly. However, we observe the simple fact that if $HH^*(A)$ is purely even, then its Hochschild homology $HH_*(A)$ is either purely even, or purely odd, since $A$ is Calabi-Yau. This implies that $HH_*(A)\otimes_\mathbb{K} R/\mathfrak{m}^k$ can not support non-zero differential just by degree reason. The statement for the negative cyclic homology is proved in the same way.
\end{Proof}

\begin{Theorem}~\label{thm:nc-vshs}
Let $A$ be a $\Z/2\Z$-graded compact, smooth, cyclic $A_\infty$ algebra, such that its Hochschild cohomology $HH^*(A)$ is concentrated at purely even degree. Let $\cA$ be its versal formal deformation over a regular local ring $(R,\mathfrak{m})$. Then the triple $\big(HC^-_*(\cA), \nabla, \langle-,-\rangle_{\hres}\big)$ forms a polarized VSHS.
\end{Theorem}

\begin{Proof}
The definition of $HC^-_*(\cA)$ using the $\mathfrak{m}$-adic completed negative cyclic homology complex allows the use of spectral sequence argument to deduce that
\[ HC^-_*(\cA) \cong \varprojlim_{k} HC^-_*(\cA/\mathfrak{m}^k\cA).\]
The locally freeness of $HC^-_*(\cA)$ follows then from the previous Lemma. The non-degeneracy of the pairing follows from that of the central fiber by Nakayama's Lemma.
\end{Proof}




\section{Kontsevich's deformation quantization and cyclic minimal models}~\label{sec:kontsevich}

Let $W\in \mathbb{K}[[x_1,\cdots,x_N]]$ be a regular function with an isolated singularity at origin.  A result of Dyckerhoff~\cite{Dyc} shows that the category of matrix factorizations $\MF(W)$ is compactly generated by the Koszul factorization $\mathbb{K}^{{\sf stab}}$. Furthermore, it was also proved in {\sl loc. cit.}  that $\MF(W)$ admits a Calabi-Yau structure. Kontsevich-Soibelman~\cite{KonSoi} showed that this implies the exitence of a cyclic minimal model $A_\infty$ algebra quasi-isomorphic to ${{\sf End}}(\mathbb{K}^{{\sf stab}})$. In this section, we apply Kontsevich's formula in deformation quantization to write down an explicit cyclic $A_\infty$ minimal model of $\End(\mathbb{K}^{{\sf stab}})$. Furthermore, if $G\subset {{\sf SL_N}}$ be a finite group of linear symmetries of $W$ preserving the top form $dx_1\wedge\cdots\wedge dx_N$, this cyclic $A_\infty$ minimal model is automatically $G$-equivariant. Thus, we may apply the smash product construction to obtain a cyclic $A_\infty$ model for the category $\MF(W,G)$ of $G$-equivariant matrix factorizations.

\paragraph{{\bf Deformation quantization in odd variables.}}
Let $A:=\mathbb{K}[ \epsilon_1,\ldots,\epsilon_N ] $ be the super-symmetric algebra generated by $N$ odd elements. Denote by $\delta: A[1]^{\otimes 2} \ra A[1]$ the shifted product map. There are two differential graded Lie algebras (DGLA) naturally associated to $A$: 
\begin{itemize}
\item[(1.)] the shifted Hochschild cochain complex of $A$ endowed with Hochschild differential $[\delta,-]$ and the Gerstenhaber bracket $\big( CC^*(A)[1],[\delta,-], [-,-]_G \big)$
\item[(2.)] the space of poly-vector fields of $A$ (see Equation~\ref{eq:poly-vector}) endowed with the trivial differential and the Schouten bracket $\big( T_{{\sf poly}}(A)[1],0, [-,-]_S\big)$. 
\end{itemize}
A fundamental result due to Kontsevich~\cite{Kon2} is the following formality

\begin{Theorem}
There exists an $L_\infty$ quasi-isomorphism
\[ \mathcal{U}: \big( T_{{\sf poly}}(A)[1],0, [-,-]_S\big) \ra \big( CC^*(A)[1], [\delta,-],  [-,-]_G \big).\]
\end{Theorem}

\medskip
\noindent Let $W(x_1,\cdots, x_N) \in \mathbb{K}[[x_1,\cdots,x_N]]$ be a formal power series such that there is no constant and linear terms. Replacing $x_j$ by $\partial_j$ we obtain a polyvector field
\[ W(\partial_1,\cdots,\partial_N) \in T_{{\sf poly}}(A)[1]\]
with constant coefficients. But any such polyvector fields trivially satisfy the Maurer-Cartan equation of the DGLA  $\big( T_{{\sf poly}}(A)[1],0, [-,-]_S\big)$, i.e. the following equation holds:
\[ [W, W]_S=0.\]
Thus, its push-forward
\[ \mathcal{U}_* W := \sum_{k\geq 1} \frac{1}{k!} \mathcal{U}_k(W,\cdots, W)\]
is a Maurer-Cartan element of $\big( CC^*(A)[1],[\delta,-], [-,-] \big)$~\footnote{Note that this infinite sum is well-defined.}. Such a Maurer-Cartan element defines an $A_\infty$ structure on $A$, corresponding to the coderivation $\delta+\mathcal{U}_* W$. More explicitly, unwinding Kontsevich's formula of $\mathcal{U}_k$ yields the higher products $\mu_n: A[1]^{\otimes n} \ra A[1]$ are given by
\[ \mu_n(a_1,\cdots, a_n)= \sum_{k\geq 1} \frac{1}{k!}\cdot \sum_{\Gamma\in G(k,n)} w_\Gamma \mathcal{U}_\Gamma(a_1,\cdots, a_n).\]
This formula needs several explanations. First, the set $G(k,n)$ is the so-called admissible graphs. Recall from~\cite{Kon2} an admissible graph $\Gamma$ of type $(k,n)$ is an oriented graph with a set of vertices $V_\Gamma:=\left\{ 1,\ldots,k\right\} \sqcup \left\{ \overline{1},\ldots,\overline{n}\right\}$, and a set of edges $E_\Gamma\subset V_\Gamma\times V_\Gamma$.  Vertices in $\left\{ 1,\ldots,k\right\}$ are called the first type, and vertices from $\left\{ \overline{1},\ldots,\overline{n}\right\}$ are called the second type. We require that every edge $(v_1,v_2)$ in $\Gamma$ must starts at a vertex $v_1$ of the first type, and ends at a vertex of the second type~\footnote{In general, there is no restriction on the ending vertex of an edge. But in our case, since $W$ is viewed as a differential operator with constant coefficients, if $\Gamma$ has an interior edge, its contribution gives zero.}.

Secondly, the weight $w_\Gamma$ of an admissible graph is given by formula
\[ w_\Gamma:=\frac{1}{(2\pi)^{2k+n-2}}\cdot \int_{C_{k,n}} \bigwedge_{e\in E_\Gamma} d\phi_e.\]
Here $C_{k,n}$ is the moduli space of $k+n$ points in the upper half plane $\mathcal{H}$ with $k$ interior points labeled by $\left\{ 1,\ldots,k\right\} $  and $n$ points on the real line labeled by $\left\{ \overline{1},\ldots,\overline{n}\right\}$ from left to right. The function $\phi_e$ is the Kontsevich's angle function associated with the edge $e$. It is given by $2$ times the angle formed by the edge $e$ with the positive real axis. The space $C_{k,n}$ is naturally an oriented manifold of dimension $2k+n-2$.  A typical configuration $\Gamma\in C_{2,4}$ is illustrated in the following picture.
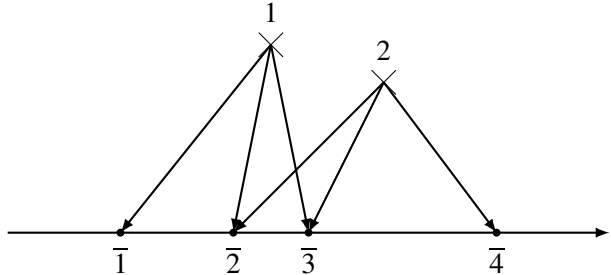
\[\begin{tikzpicture}[baseline={(current bounding box.center)},scale=0.5]
\draw (-1,5) node[cross=5pt,label=above:{1}] {};
\draw (2,4) node[cross=5pt,label=above:{2}] {};
\node[circle,fill=black,inner sep=0pt,minimum size=3pt,label=below:{$\overline{1}$}] at (-5,0) {};
\node[circle,fill=black,inner sep=0pt,minimum size=3pt,label=below:{$\overline{2}$}] at (-2,0) {};
\node[circle,fill=black,inner sep=0pt,minimum size=3pt,label=below:{$\overline{3}$}] at (0,0) {};
\node[circle,fill=black,inner sep=0pt,minimum size=3pt,label=below:{$\overline{4}$}] at (5,0) {};
\draw [thick,->,>=latex,] (-1,5) to (-5,0);
\draw [thick,->,>=latex] (-1,5) to (-2,0);
\draw [thick,->,>=latex] (-1,5) to (0,0);
\draw [thick,->,>=latex] (2,4) to (-2,0);
\draw [thick,->,>=latex] (2,4) to (0,0);
\draw [thick,->,>=latex] (2,4) to (5,0);
\draw [thick,->,>=latex] (-8,0) to (8,0);
\end{tikzpicture}\]
The integral above depends on a choice of ordering of the set $E_\Gamma$. We may choose the following canonical ordering. For the $i$-th interior marked point, we set ${{\sf Star}}(i)$ to be the set of out-going (which is the same as all in our case) edges at the vertex $i$ in the Graph $\Gamma$. The set ${{\sf Star}}(i)$ is ordered canonically according to the ordering of the ending points of the edges on the real line. A canonical ordering of $E_\Gamma$ is then given by
\[ {{\sf Star}}(1),\cdots,{{\sf Star}}(k).\]
With the canonical ordering, the weight $w_\Gamma$ may be computed explicitly. The answer turned out to be rational numbers.
\begin{Lemma}
The weight of an admissible graph is given by
\[ w_\Gamma = \frac{1}{|\left\{\Gamma-{{\sf Shuffles}}\right\}|}\cdot \prod_{i=1}^k \frac{1}{|{{\sf Star(i)}}|!}.\]
An element $\sigma\in S_n$ is called a $\Gamma$-Shuffle if for each $1\leq i\leq k$ the restriction $\sigma|_{{{\sf Star(i)}}}: {{\sf Star(i)}} \rightarrow \left\{ 1,2,\cdots,n\right\}$ is an order-preserving map. Here we identified an edge $e\in {{\sf Star(i)}}$ with the its end point in $\left\{ 1,2,\cdots,n\right\}$.
\end{Lemma}

\begin{Proof}
Let $\Gamma\in G(k,n)$ be an admissible graph. It determines an embedding 
\[ \iota_\Gamma: C_{k,n}\hookrightarrow \Delta^{|{{\sf Star(1)}}|}\times\cdots\times \Delta^{|{{\sf Star(k)}}|}.\]
Here each $\Delta^{|{{\sf Star(i)}}|}$ is the standard simplex of dimension $|{{\sf Star(i)}}|$. The weight of $\Gamma$ is then just given by $\int_{\iota_\Gamma(C_{k,n})} d\mu$ with $d\mu$ the Euclidean measure.

Consider an admissible graph as a map $\Gamma:  \left\{ \overline{1},\cdots,\overline{k}\right\} \ra  \left\{ 1,2,\cdots,n\right\}$. And let $\sigma$ be a $\Gamma$-Shuffle. Then $\sigma\Gamma$ is another admissible graph.  The point is that we have a decomposition 
\[ \coprod_\sigma \iota_{\sigma\Gamma}(C_{k,n}) = \Delta^{|{{\sf Star(1)}}|}\times\cdots\times \Delta^{|{{\sf Star(k)}}|}.\]
Furthermore, each $\sigma\in S_n$ acts on $\Delta^{|{{\sf Star(1)}}|}\times\cdots\times \Delta^{|{{\sf Star(k)}}|}$ by permutation of the coordinates, which identifies the weights $w_\Gamma=w_{\sigma\Gamma}$. Since each simplex $\Delta^{|{{\sf Star(i)}}|}$ has volume equal to $\frac{1}{|{{\sf Star(i)}}|!}$, we obtain the weight formula in the Lemma. 
\end{Proof}

\medskip
The action $\mathcal{U}_\Gamma$ of an admissible graph is defined as follows. For each monomial $\partial_I=\partial_{i_1}\cdots\partial_{i_l}\in T_{{\sf poly}}(A)$, we set
\[ \widetilde{\partial_I}:= \sum_{\sigma\in S_l} \partial_{i_{\sigma(1)}}\otimes \cdots \otimes \partial_{i_{\sigma(l)}}\in {{\sf Hom}} (A, A)^{\otimes l}.\]
Extending $\mathbb{K}$-linearly this operation yields $\widetilde{W}\in \prod_l {{\sf Hom}} (A[1], A[1])^{\otimes l}$. Denote by $\widetilde{W}_l$ the $l$-th component. For each interior vertex $1\leq i\leq k$, define an operator 
\[  \widetilde{W}_{{\sf Star(i)}}: A[1]^{\otimes n} \rightarrow A[1]^{\otimes n},\]
which acts on the component $\otimes_{j, (i,\overline{j})\in {{\sf Star(i)}}} A[1]$ by $\widetilde{W}_{|{\sf Star(i)}|}$, and by the identity map on the rest tensor components. The action $\mathcal{U}_\Gamma$ is given by
\[ \mathcal{U}_\Gamma (a_1,\cdots,a_n):= M\big( \widetilde{W}_{{{\sf Star(1)}}} \circ\cdots\circ \widetilde{W}_{{{\sf Star(k)}}}(a_1,\cdots,a_n)\big).\]
Here $M$ is the product operator of multiplying $n$ elements in the exterior algebra $A[1]=\mathbb{K}[\epsilon_1,\cdots,\epsilon_N]$. Explicitly, it is defined by
\[ M(a_1,\cdots,a_n)= (-1)^{\frac{n(n-1)}{2}+(n-1)|a_1|+(n-2)|a_2|+\cdots+|a_{n-1}|} a_1\wedge\cdots\wedge a_n.\]
\paragraph{{\bf Explicit formulas.}} Putting all the above discussions together we obtain the following expression of the $A_\infty$ multiplication of $A^W$:
\begin{align*}
 \mu_n(a_1,\cdots,a_n) &=\sum_{k\geq 1} \frac{1}{k!}\cdot \sum_{\Gamma\in G(k,n)}  w_\Gamma\cdot M\big( \widetilde{W}_{{{\sf Star(1)}}} \circ\cdots\circ \widetilde{W}_{{{\sf Star(k)}}}(a_1,\cdots,a_n)\big)\\
 w_\Gamma &= \frac{1}{|\left\{\Gamma-{{\sf Shuffles}}\right\}|}\cdot \prod_{i=1}^k \frac{1}{|{{\sf Star(i)}}|!}
 \end{align*}
Let us illustrate the above formula in the case when $a_1=\epsilon_{j_1},\cdots,a_n=\epsilon_{j_n}$ are linear tensors.  There is a  unique graph that can contribute:
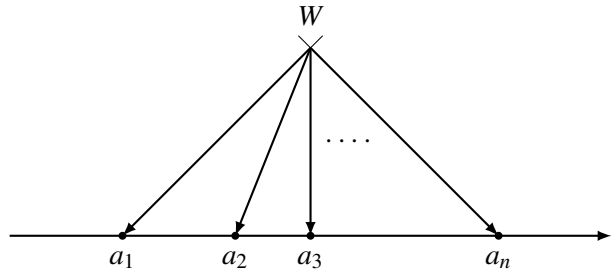
\[\begin{tikzpicture}[baseline={(current bounding box.center)},scale=0.5]
\draw (0,5) node[cross=5pt,label=above:{$W$}] {};
\node[circle,fill=black,inner sep=0pt,minimum size=3pt,label=below:{$a_1$}] at (-5,0) {};
\node[circle,fill=black,inner sep=0pt,minimum size=3pt,label=below:{$a_2$}] at (-2,0) {};
\node[circle,fill=black,inner sep=0pt,minimum size=3pt,label=below:{$a_3$}] at (0,0) {};
\node[circle,fill=black,inner sep=0pt,minimum size=3pt,label=below:{$a_n$}] at (5,0) {};
\node at (1,2.5) {$\cdots\cdot$};
\draw [thick,->,>=latex] (0,5) to (-5,0);
\draw [thick,->,>=latex] (0,5) to (-2,0);
\draw [thick,->,>=latex] (0,5) to (0,0);
\draw [thick,->,>=latex] (0,5) to (5,0);
\draw [thick,->,>=latex] (-8,0) to (8,0);
\end{tikzpicture}\]
Its weight is $\frac{1}{n!}$. Assuming that for each $1\leq j\leq N$, $\epsilon_j$ appears $k_j$ times in $a_1,\cdots, a_n$, then there are $k_1!\cdots k_N!$ non-zero contributions in $\widetilde{W}_n$. This gives the following formula:
 \[ \mu_n(\epsilon_{j_1},\cdots,\epsilon_{j_n})=\frac{1}{n!}\partial_{x_{j_1}}\cdots\partial_{x_{j_n}} W |_{x=0} \cdot \bone.\]
This refines a formula of Dyckerhoff which gives symmetrizations of the above formula~\cite{Dyc}.

\paragraph{{\bf Comparison of $A^W$ with the minimal model of $\End(\mathbb{K}^{{\sf stab}})$.}} To avoid long notations, we use the symbol $A^W$ to mean the $A_\infty$ algebra defined by $\delta+\mathcal{U}_*W$. To relate to the category of matrix factorizations, we shall prove that the $A_\infty$ algebra $A^W$ is a minimal model of $\End(\mathbb{K}^{{\sf stab}})$. 

Let $\mathfrak{g}$ be a differential graded Lie algebra. For each negative integer $n\in \left\{ -1,-2,\cdots\right\}$, we set the smart truncation 
\[ \tau^{\leq n}\mathfrak{g}:= \cdots\ra \mathfrak{g}^{n-2}\ra \mathfrak{g}^{n-1} \ra \ker d_n \ra 0, \]
where $d_n: \mathfrak{g}^n \ra \mathfrak{g}^{n+1}$ is the differential at the $n$-th degree. Observe that the truncation $\tau^{\leq n}\mathfrak{g}$ is a differential graded Lie subalgebra of $\mathfrak{g}$, for any negative integer $n$. Kontsevich's  $L_\infty$ homomorphism $$\mathcal{U}: \big( T_{{\sf poly}}(A)[1],0, [-,-]_S\big) \ra \big( CC^*(A)[1], [\delta,-], [-,-]_G\big)$$ restricts to a $L_\infty$ homomorphism
\[ \tau^{\leq -1} \mathcal{U}: \tau^{\leq -1}T_{{\sf poly}}(A)[1]\mathfrak{g} \ra \tau^{\leq -1} CC^*(A)[1].\]
This is well-defined since $\mathcal{U}_k$ is of degree $1-k$ with $k\geq 1$, and $\mathcal{U}_1$ is a map of complexes. 

The truncation $\tau^{\leq -1}  T_{{\sf poly}}(A)[1] = \mathbb{K}[[\partial_1,\cdots,\partial_N]]$ is given by polyvector fields with constant coefficients concentrated at degree $-1$. It is endowed with trivial differential and trivial Lie bracket. The other truncation $\tau^{\leq -1} CC^*(A)[1]$ has a {{\sl pro-nilpotent} DGLA structure, given by the decreasing filtration $\tau^{\leq -1} CC^*(A)[1]\supset \tau^{\leq -2} CC^*(A)[1] \supset \tau^{\leq -3} CC^*(A)[1] \supset \cdots$.

Let $\mathfrak{h}$ be a negatively graded pro-nilpotent DGLA. We define the set of negative Maurer-Cartan elements by
\[ \widetilde{\MC^-}(\mathfrak{h}):=\left\{ \alpha=\alpha_{-1}+\alpha_{-3}+\alpha_{-5}+\cdots \in \prod_{n\; \mbox{odd}} \mathfrak{h}_n \mid d\alpha+\frac{1}{2} [\alpha,\alpha]=0.\right\}\]
Two negative Maurer-Cartan elements $\alpha_1$ and $\alpha_2$ are said to be gauge equivalent if there exists an element 
\begin{align*}
h & =h_{-2}+h_{-4}+\cdots\in \prod_{n\; \mbox{even}} \mathfrak{h}_n,\; \mbox{ such that}\\
\alpha_2 &=\exp(h)* \alpha_1 := \exp({{\sf ad}}(h)) \alpha_1 - \frac{\exp({{\sf ad}}(h))-1}{{{\sf ad}}(h)} dh.
\end{align*}
We set $\MC^-(\mathfrak{h}):= \widetilde{\MC^-}(\mathfrak{h}) / \mbox{gauge equivalences}$.

\begin{Lemma}
The $L_\infty$ quasi-isomorphism 
\[\tau^{\leq -1}\mathcal{U}: \mathbb{K}[[\partial_1,\cdots,\partial_N]] \ra \tau^{\leq -1} CC^*(A)[1]\] 
induces an bijection of sets
\[ \big(\tau^{\leq -1}\mathcal{U}\big)_*: \mathbb{K}[[\partial_1,\cdots,\partial_N]] \ra \MC^-\big( \tau^{\leq -1} CC^*(A)[1]\big).\]
\end{Lemma}

\begin{proof}
It is easy to see that the map $\big(\tau^{\leq -1}\mathcal{U}\big)_*$ is injective. Indeed, the equation 
\[ \big(\tau^{\leq -1}\mathcal{U}\big)_*W_1=\big(\tau^{\leq -1}\mathcal{U}\big)_*W_2\]
implies that $\mathcal{U}_1(W_1)$ is cohomologous to $\mathcal{U}_1(W_2)$. But $\mathcal{U}_1$ is the HKR map which is a quasi-isomorphism. Hence the equality $[\mathcal{U}_1(W_1)]=[\mathcal{U}_1(W_2)]$ implies that $W_1=W_2$. 

For the surjectivity, it suffices to show that two Maurer-Cartan elements $a=a_{-1}+a_{-3}+\cdots$ and $b=b_{-1}+b_{-3}+\cdots$ such that $a_{-1}$ is cohomologous to $b_{-1}$ are gauge equivalent. The idea is to use the degree filtration to inductively construct $h=h_{-2}+h_{-4}+\cdots$.  
Indeed, since $a_{-1}$ and $b_{-1}$ are cohomologous, we can find $h_{-2}$ such that
\[ dh_{-2}=a_{-1}-b_{-1}.\]
Now assuming that we have constructed $h_k$ for all $k\geq n$ for some negative even integer $n$, we would like to construct $h_{n-2}$. Let us write $h_{\geq n}:= h_{-2}+\cdots+h_{n}$. By assumption Maurer-Cartan element $b':=\exp(h_{\geq n})*a$ agrees with $b$ at degrees that are greater or equal to $n+1$. But for any Maurer-Cartan element $x$ we have
\[dx_{n-1}=-\frac{1}{2}\sum_{i+j=n}[x_i,x_j],\]
which shows that $dx_{n-1}$ is determined by elements of strictly bigger degrees. Hence we conclude that
\[ db'_{n-1}=db_{n-1}.\]
Thus we deduce that $d(b'_{n-1}-b_{n-1})=0$. Since the cohomology of $\tau^{\leq -1} CC^*(A)[1]$ vanishes at degrees that are less or equal to $-2$, there exists $h_{n-2}$ such that
\[ dh_{n-2}=b'_{n-1}-b_{n-1}.\]
This implies the Maurer-cartan element $\exp(h_{\geq n}+h_{n-2})*a$ agrees with $b$ in degrees that are greater or equal to $n-1$.
\end{proof}

The lemma easily implies the following

\begin{Corollary}~\label{cor:minimal}
Let $[\alpha^{{\sf LG}}]\in \MC\big( \tau^{\leq -1} CC^*(A)[1]\big)$ be the point in the Maurer-Cartan moduli space of $\big(\tau^{\leq -1} CC^*(A)[1], [\delta,-], [-,-]\big)$ corresponding to the minimal model $A_\infty$ algebra of $\End(\mathbb{K}^{{\sf stab}})$ defined by Dyckerhoff. Then we have
\[ [\alpha^{{\sf LG}}] = [\tau^{\leq -1}\mathcal{U}_*W].\]
That is, the $A_\infty$ structure we obtained from Kontsevich's deformation quantization is gauge equivalent to the minimal model $A_\infty$ structure of $\End(\mathbb{K}^{{\sf stab}})$ obtained from homological perturbation technique.
\end{Corollary}

\begin{proof}
By the lemma above, it suffices to prove that the degree $-1$ part of the $A_\infty$ structure $\alpha^{{\sf LG}}$ is cohomologous to $\mathcal{U}_1(W)$. For this, we may use the splitting map of the HKR-isomorphism given by
\[ \phi_n \mapsto \phi_n(\epsilon_{j_1},\cdots,\epsilon_{j_n})\partial_{j_1}\cdots\partial_{j_n}.\]
For $\alpha^{{\sf LG}}_{-1}$, this splitting map gives $W$ by Dyckerhoff's formula~\cite{Dyc}. For $\mathcal{U}_1(W)$ this also gives $W$ since $\mathcal{U}_1$ is the HKR map.
\end{proof}

\paragraph{{\bf Hochschild cohomology of $(A,\delta+\mathcal{U}_*W)$.}}  A particularly nice feature of Kontsevich's explicit formula is that it preserves the cup product structure on the tangent cohomology. Namely, it was observed by Kontsevich that the tangent map at $W$ defines an isomorphism of algebras 
\begin{align*}
 d\mathcal{U}_W &:H^*\big( T_{{\sf poly}}(A), [W,-]\big) \ra HH^*(A^W),\\
 d\mathcal{U}_W ([h]) &:=\sum_{k\geq 0} \frac{1}{k!} \mathcal{U}_{k+1}(W,\cdots,W,h).
     \end{align*}
Note that this result does not require that $W$ have isolated singularities. If $W$ has an isolated singularity at origin, the above isomorphism reduces to the well-known fact that the Hochschild cohomology of $A$ is isomorphic to the Jacobi algebra ${{\sf Jac}}(W)$. Interestingly, even in this case, the above algebra map is {{\sl not}} given by the Hochschild-Konstant-Rosenberg isomorphism which is only the first term of $d\mathcal{U}_W$. It would be very interesting to explore the Duflo type correction terms in this setting.

\paragraph{{\bf Cyclicity of $A^W$.}} Another nice feature of Kontsevich's formula is that the minimal model $A_\infty$ structure $\delta+\mathcal{U}_*W$ is automatically cyclic, with respect to the trace map
\[ \Tr: \mathbb{K}[ \epsilon_1,\ldots,\epsilon_N ]  \ra \mathbb{K}, \;\;\; \Tr(\epsilon_I)=\begin{cases} 1, \mbox{\;\;\; if $I=(1,2,3,\cdots,N)$}\\ 0, \mbox{\;\;\; otherwise.} \end{cases}\]
Indeed, the trace map induces a bilinear form on $\mathbb{K}[ \epsilon_1,\ldots,\epsilon_N ]$ by sending $a\otimes b \mapsto \Tr(a\wedge b)$. The non-degeneracy of this bilinear form further induces an isomorphism $\iota: A\ra A^\vee$. We use $\iota$ to get an isomorphism
\[ CC^*(A) \cong \big( CC_*(A)\big)^\vee.\]
The dual of the Connes differential $B^\vee$ pulls back to $CC^*(A)$ which we denote by ${{\sf div}}_{\Tr}$, called the divergence operator associated to the trace map $\Tr$.
Since $W$ has constant coefficients, we have ${{\sf div}}_{\Tr}(W)=0$. It was observed by Felder-Shoikhet~\cite{FelSho} (see also Willwachera-Calaque~\cite{WilCal} for the non-divergence free case) that for divergence free polyvector fields, the push-forward under Kontsevich's $L_\infty$ homomorphism $\mathcal{U}$ is cyclic, in the sense that
\[ \Tr\big( \mu_n(a_1,\cdots, a_n)\wedge a_{n+1}\big) = (-1)^{|a_{n+1}|'(|a_1|'+\cdots+|a_n|')} \Tr\big( \mu_n(a_{n+1},a_1,\cdots,a_{n-1})\wedge a_n\big).\]
This is precisely the cyclicity of $\mu$ with respect to the trace map $\Tr$.

\paragraph{{\bf Deformations of $A^W$.}} Assume that $W$ has an isolated singularity at origin. Consider a deformation of $W$ depending on parameters $t_1,\cdots,t_\mu$ defined by $\mathscr{W}=W+t_1\phi_1+\cdots+t_\mu\phi_\mu$ with $(\phi_1,\cdots,\phi_\mu)$ a basis of the Jacobian ring ${{\sf Jac}}(W)$. Extending Kontsevich's formality map $t$-linearly yields a morphism
\[ \mathcal{U}: T_{{\sf poly}}(A)[[t_1,\cdots,t_\mu]][1] \ra C^*(A)[[t_1,\cdots,t_\mu]][1].\]
The push-forward $\mathcal{U}_*(\mathscr{W})$ yields a flat family of $A_\infty$ structure $A^{{\mathscr{W}}}$ over the commutative ring $R=\mathbb{K}[[t_1,\cdots,t_\mu]]$. Observe that by construction, this family $A^{{\mathscr{W}}}$ is a versal deformation of $A^W$. Theorem~\ref{thm:nc-vshs} applied to the family $A^{{\mathscr{W}}}$ yields a VSHS:
\[ \cV^{\MF(W)}:=\Big( HC^-_*(A^{\mathscr{W}}), \nabla, \langle-,-\rangle_{{\sf hres}}\Big).\]

\paragraph{{\bf Symmetries in $G\subset \SL_N$.}} Let $G\subset \SL_N$ be a finite subgroup of linear symmetries of $W$. Yet another nice property of Kontsevich's formula of $\mathcal{U}$ is that it is equivariant under all affine transformations of $\C^N$. In particular, $\mathcal{U}$ is $G$-equivariant, which implies that the higher products $\mu_n$ are $G$-equivariant maps. Furthermore, since $G\subset \SL_N(\C)$, it also preserves the trace map, which implies that $G$ is a group of symmetries of the cyclic $A_\infty$ algebra $(A^W,\Tr)$.

This enables us to form another cyclic $A_\infty$ algebra $A^W\rtimes G$, the semi-direct product algebra. As a vector space $A^W\rtimes G= A^W \otimes_\mathbb{K} \mathbb{K}[G]$, its higher products are defined by
\[ \mu_n( a_1\otimes g_1,\cdots, a_n\otimes g_n):= \mu_n\big(a_1,g_1(a_2),g_1g_2(a_3),\cdots,g_1g_2\ldots g_{n-1}(a_n)\big)\otimes g_1\ldots g_n.\]
The trace map is defined by
\[ \Tr(\epsilon_I\otimes g) = \begin{cases} 1, \mbox{\;\;\; if $I=(1,2,3,\cdots,N)$ and $g=\id_G$ }\\ 0, \mbox{\;\;\; otherwise.} \end{cases}\]
The semi-direct product algebra naturally arises from the category $\MF(W,G)$ of $G$-equivariant matrix factorizations of $W$. Indeed, Polishchuk-Vaintrob~\cite{PolVai} proved that $\mathbb{K}^{{\sf stab}}\otimes_\mathbb{K}\mathbb{K}[G]$ compactly generates the category $\MF(W,G)$. Its endomorphism differential graded algebra is 
\[ \End_{\MF(W,G)}\big( \mathbb{K}^{{\sf stab}}\otimes_\mathbb{K}\mathbb{K}[G]\big) = \End_{\MF(W)}(\mathbb{K}^{{\sf stab}}) \rtimes G.\]
Passing to the minimal model yields precisely $A^W\rtimes G$.

\section{Tsygan formality and Comparison of VSHS's}~\label{sec:tsygan}

In this section, we prove the main result of the paper that the deformed Tsygan formality map is an isomorphism of VSHS's between $\cV^{{\MF(W)}}$ and $\cV^{{\sf Saito}}$.

\paragraph{{\bf Tsygan formality.}} Recall our notation that $A=\mathbb{K}[\epsilon_1,\cdots,\epsilon_N]$, the super-commutative algebra generated by $N$ odd variables in cohomological degree $1$. And there is an $L_\infty$ quasi-isomorphism
\[ {{\mathcal{U}}}: T_{{\sf poly}} (A) [1] \rightarrow CC^*(A)[1],\]
from the DGLA of poly-vector fields to the Hochschild cochain complex. Let $\phi\in CC^*(A)[1]$ be a homogeneous Hochschild cochain. The assignment $\phi\mapsto L_\phi$ endows $CC_*(A)$ a differential graded Lie module structure over $CC^*(A)[1]$. Precomposing with the $L_\infty$ quasi-isomorphism $\mathcal{U}$ gives $CC_*(A)$ an $L_\infty$-module structure over $T_{{\sf poly}} (A)$. On the commutative side, the space $\Omega^*_A$ of differential forms is naturally a graded Lie-module over $T_{{\sf poly}} (A)[1]$ through the standard Lie derivative action. Tsygan formality asserts that there exists an $L_\infty$ module quasi-isomorphism
\[ \mathcal{U}^{Sh}: CC_*(A) \rightarrow \Omega^*_A.\]
Following Willwacher~\cite{Wil}, here the superscript ``Sh" is for Shoikhet who wrote down this $L_\infty$-module map~\cite{Sho}. The Shoikhet quasi-isomorphism is defined in the same fashion as the map $\mathcal{U}$ in the case of Kontsevich formality. Indeed, for a Hochschild chain $a_0|\cdots|a_n$, and Hochschild cochains $\phi_1,\cdots,\phi_k\in CC^*(A)[1]$, the differential form $\mathcal{U}_k^{Sh} (\phi_1,\cdots,\phi_k; a_0|\cdots |a_n)$ is given by
\[  \mathcal{U}_k^{Sh} (\phi_1,\cdots,\phi_k; a_0|\cdots |a_n):= 
\sum_{\substack{\Gamma\in G'(k+1,n+1)\\ j\geq 0\\ I=(i_1,\cdots,i_j)\in \left\{ 1,\cdots, N\right\}^j}} w_\Gamma \mathcal{U}_\Gamma(\partial_I,\phi_1,\cdots,\phi_k;a_0|\cdots|a_n) d\epsilon_I.\]
In this formula, for each multi-index $ I=(i_1,\cdots,i_j)$, we used the abbreviations $\partial_I:=\partial_{i_1}\cdots\partial_{i_j}$ and $d\epsilon_I=d\epsilon _{i_1}\cdots d\epsilon_{i_j}$. Several explanations of this formula are in order. 
\begin{itemize}
\item[--] First, the weight $w_\Gamma$ is defined by the following integral
\[ w_\Gamma:=\int_{D_{k+1,n+1}} \bigwedge_{e\in E_\Gamma} d\theta_e,\]
with $\theta_e$ the Shoikhet angle function associated with the edge $e$, while $D_{k+1,n+1}$ is the configuration space of $k+1$ interior points and $n+1$ boundary points on the unit disk. The $n+1$ boundary marked points are labled by $\left\{ \overline{0}, \cdots,\overline{n}\right\}$ in the counter-clockwise orientation, and the $k+1$ interior marked points are labeld by $\left\{ 0,\cdots,k\right\}$. We may fix the zeroth boundary marked point at $(1,0)$, and the zeroth interior marked point at $(0,0)$. 
\item[--] the set $G'(k+1,n+1)$ consists of admissible directed graphs, such that a graph $\Gamma\in G'(k+1,n+1)$ is admissible if it satisfies the conditions as in Kontsevich's case, plus that the first interior point $0$ can not be the ending vertex of an edge. 
\item[--]The action $\cU_\Gamma$ is defined by the same rule as in Kontsevich's case, by putting the poly-vector field $\partial_I$ at the $0$-th interior marked point, and $a_j$ at the $\overline{j}$-th boundary marked point for $0\leq j\leq n$. The following picture illustrates the definition of $\cU_{\Gamma}(\partial_I,\phi_1,\cdots,\phi_k;a_0|\cdots|a_n)$. 
\end{itemize}
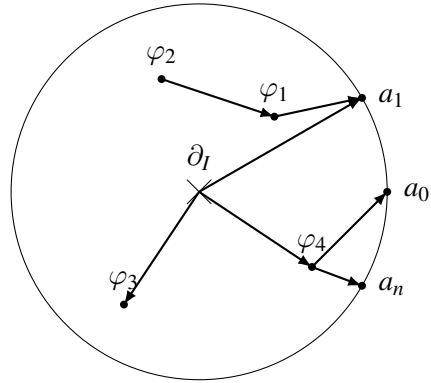
\[ \begin{tikzpicture}[baseline={(current bounding box.center)},scale=0.5]
\draw (0,0) circle (5);
\node[circle,fill=black,inner sep=0pt,minimum size=3pt,label=right:{$a_0$}] at (5,0) {};
\node[circle,fill=black,inner sep=0pt,minimum size=3pt,label=right:{$a_1$}] at (4.33,2.5) {};
\node[circle,fill=black,inner sep=0pt,minimum size=3pt,label=right:{$a_n$}] at (4.33,-2.5) {};
\draw (0,0) node[cross=5pt,label=above:{$\partial_I$}] {};
\node[circle,fill=black,inner sep=0pt,minimum size=3pt,label=above:{$\phi_1$}] at (2,2) {};
\node[circle,fill=black,inner sep=0pt,minimum size=3pt,label=above:{$\phi_2$}] at (-1,3) {};
\node[circle,fill=black,inner sep=0pt,minimum size=3pt,label=above:{$\phi_3$}] at (-2,-3) {};
\node[circle,fill=black,inner sep=0pt,minimum size=3pt,label=above:{$\phi_4$}] at (3,-2) {};
\draw [thick,->,>=latex] (0,0) to (4.33,2.5);
\draw [thick,->,>=latex] (0,0) to (-2,-3);
\draw [thick,->,>=latex] (0,0) to (3,-2);
\draw [thick,->,>=latex] (3,-2) to (5,0);
\draw [thick,->,>=latex] (3,-2) to (4.33,-2.5);
\draw [thick,->,>=latex] (-1,3) to (2,2);
\draw [thick,->,>=latex] (2,2) to (4.33,2.5);
\end{tikzpicture}\]

\paragraph{{\bf Deformations of the Tsygan formality map.}} Let $W\in T_{{\sf poly}}(A)$ be a poly-vector field with constant coefficients, i.e. corresponding to a function $W\in \mathbb{K}[[ x_1,\cdots,x_N]]$. Since such a $W$ is a Maurer-Cartan element of $T_{{\sf poly}}(A)[1]$, it induces a deformation of the Tsygan formality map
\begin{align*}
 {\mathcal{U}}^{Sh,W} &: CC_*(A^W) \rightarrow \big( \Omega^*_A, L_W\big).
 \end{align*}
 We shall only be intersted with the first deformed structure map $\mathcal{U}^{Sh,W}_0$ explicitly given by
 \[ \mathcal{U}^{Sh,W}_0(-)=\sum_{k\geq 0} \frac{1}{k!} \mathcal{U}^{Sh}_k(\underbrace{W,\cdots,W}_{\mbox{$k$-copies}};-).\]
That is, we insert $W$ at all interior marked points except the zeroth one. Its $k$-th term is illustrated in the following picture where we insert $k$-copies of $W$ in the interior of the disc.
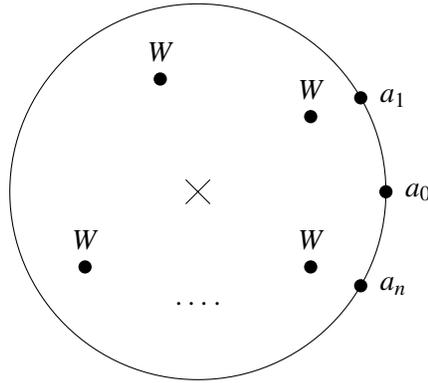
\[ \begin{tikzpicture}[baseline={(current bounding box.center)},scale=0.5]
\draw (0,0) circle (5);
\node[circle,fill=black,inner sep=0pt,minimum size=5pt,label=right:{$a_0$}] at (5,0) {};
\node[circle,fill=black,inner sep=0pt,minimum size=5pt,label=right:{$a_1$}] at (4.33,2.5) {};
\node[circle,fill=black,inner sep=0pt,minimum size=5pt,label=right:{$a_n$}] at (4.33,-2.5) {};
\draw (0,0) node[cross=5pt] {};
\node[circle,fill=black,inner sep=0pt,minimum size=5pt,label=above:{$W$}] at (3,2) {};
\node[circle,fill=black,inner sep=0pt,minimum size=5pt,label=above:{$W$}] at (-1,3) {};
\node[circle,fill=black,inner sep=0pt,minimum size=5pt,label=above:{$W$}] at (-3,-2) {};
\node[circle,fill=black,inner sep=0pt,minimum size=5pt,label=above:{$W$}] at (3,-2) {};
\node at (0,-3) {$\cdots\cdot$};
\end{tikzpicture}\]
Let us now assume that $W$ has an isolated singularity at origin. Choose a set of basis $\varphi_1,\cdots,\varphi_\mu$ for the Jacobian ring of $W$. Denote by
\[ \mathscr{W}:= W + t_1\varphi_1+\cdots+ t_\mu\varphi_\mu.\]
Replacing $W$ by $\mathscr{W}$ produces a formal deformation of the Tsygan formality map depending on the formal parameters $t_1,\cdots, t_\mu$. Denote this map by $ \mathcal{U}_0^{Sh,\mathscr{W}}$. We may extend the map $u$-linearly to obtain a map on the negative cyclic chain complex
\begin{equation}\label{eq:main} \mathcal{U}_0^{Sh,\mathscr{W}}: \big( CC_*(A^{\mathscr{W}})[[u]], b^{\mathscr{W}}+uB\big) \rightarrow \big( \Omega^*_A\otimes \mathbb{K}[[t_1,\cdots,t_\mu]][[u]], L_{\mathscr{W}}+ud_{DR}\big).\end{equation}
The fact that $\mathcal{U}_0^{Sh,\mathscr{W}}$ remains a chain map (i.e. it intertwines the Connes operator $B$ with the de Rham differential $d_{DR}$) is proved by Willwatcher~\cite[Theorem 1.3]{Wil}: we have 
\[ \mathcal{U}_0^{Sh,\mathscr{W}}B=d_{DR} \mathcal{U}_0^{Sh,\mathscr{W}}.\]
Taking cohomology of Equation~\ref{eq:main} yields an isomorphism
\[ \mathcal{U}_0^{Sh,\mathscr{W}}: HC^-_*(A^{\mathscr{W}}) \rightarrow H^*\big( \Omega^*_A\otimes \mathbb{K}[[t_1,\cdots,t_\mu]][[u]], L_{\mathscr{W}}+ud_{DR}\big).\]
The left hand side carries a VSHS as explained in Section~\ref{sec:vshs} from non-commutative geometry.

\paragraph{{\bf Koszul duality: $\Omega_A^*$ versus $\Omega_S^*$.}  For the right hand side of the above isomorphism, we shall show that it is naturally isomorphic to Saito's VSHS
$$\cV^{{\sf Saito}}=H^*\big( \Omega_S^*[[t_1,\cdots,t_\mu,u]], d\mathscr{W}+ud_{DR} \big),$$ 
assuming that $W$ has an isolated singularity at the origin. Indeed, observe the following sequence of isomorphisms
\begin{equation}~\label{eq:dualize}
 H^*\big( \Omega^*_A[[t_1,\cdots,t_\mu,u]], L_{\mathscr{W}}+ud_{DR}\big) \ra \big( \cV^{{\sf Saito}}\big)^\vee \leftarrow \cV^{{\sf Saito}}.
 \end{equation}
Here the first arrow sends
\[ \epsilon_j \mapsto (dx_j)^\vee, \;\; d\epsilon_j \mapsto (x_j)^\vee,\]
extended $\mathbb{K}[[t_1,\cdots,t_\mu]]$-linearly and $\mathbb{K}[[u]]$-sesquilinearly. The second arrow is induced by the higher residue pairing $\alpha \mapsto K^{{\mathscr{W}}}(\alpha,-)$. Note that the composition is $\mathbb{K}[[u]]$-linear as the higher residue pairing is  $\mathbb{K}[[u]]$-sesquilinearly in the second component. Since $W$ has an isolated singularity, this map induces an isomorphism modulo the ideal $(t_1,\cdots,t_\mu,u)$, which further implies that it is itself an isomorphism by Nakayama lemma. By definition, we put on the space $H^*\big( \Omega^*_A\otimes \mathbb{K}[[t_1,\cdots,t_\mu]][[u]], L_{\mathscr{W}}+ud_{DR}\big)$  the pull-back VSHS under the isomorphism~\ref{eq:dualize}. We shall need the following explicit formulas for the induced connection operators in terms of calculus of the super-commutative algebra $A$:
\begin{align*}
\nabla_{\frac{\partial}{\partial u}}& := \frac{\partial}{\partial u}+\frac{N}{2u} + \frac{\iota_{{\mathscr{W}}}}{u^2},\\
\nabla^{{\sf GM}}_{\frac{\partial}{\partial t_j}} &:= \frac{\partial}{\partial t_j} - \frac{\iota_{\phi_j}}{u}.
\end{align*}
We may now state our main result:

\begin{Theorem}
The deformed Tsygan formality map  
\[ \mathcal{U}_0^{Sh,\mathscr{W}}: HC^-_*(A^{\mathscr{W}}) \rightarrow H^*\big( \Omega^*_A\otimes \mathbb{K}[[t_1,\cdots,t_\mu]][[u]], L_{\mathscr{W}}+ud_{DR}\big).\]
is an isomorphism of VSHS's.
\end{Theorem}

\begin{Proof}
In Proposition~\ref{prop:gm}, we shall prove that $\mathcal{U}_0^{Sh,\mathscr{W}}$ intertwines with the Gauss-Manin connection $\nabla^{{\sf GM}}$ with the Getzler connection $\nabla^{{\sf Get}}$. In Proposition~\ref{prop:uconn}, we prove that $\mathcal{U}_0^{Sh,\mathscr{W}}$ also intertwines the $u$-direction connections. These two properties imply that $\mathcal{U}_0^{Sh,\mathscr{W}}$ automatically intertwines the higher resiude pairing (up to a non-zero constant in $\mathbb{K}^*$) by a uniqueness result of M. Saito~\cite{MSai}. The rest of the section is devoted to the proofs of Proposition~\ref{prop:gm} and Proposition~\ref{prop:uconn}.
\end{Proof}

\paragraph{{\bf Construction of the homotopy operators.}} Fix an integer $1\leq j\leq \mu$. We would like to prove that
\[ \nabla_{\frac{\partial}{\partial t_j}}^{{\sf GM}} \circ \mathcal{U}_0^{Sh,\mathscr{W}} - \mathcal{U}_0^{Sh,\mathscr{W}}\circ \nabla^{{\sf Get}}_{\frac{\partial}{\partial t_j}}=0,\]
after taking cohomology. To acheive this, we use a homotopy operator constructed by Cattaneo-Felder-Willwacher~\cite{CFW}. Consider the configuration space $E_{k+2,n+1}$ of $k+2$ interior marked points $(z_0,\cdots,z_{k+1})$, and $n+1$ boundary marked points $(z_{\overline{0}}, \cdots, z_{\overline{n}})$ on the unit disk, such that the first two interior points are such that $z_0=(0,0)$ and $z_1=(r,0)$ with $r\in [0,1]$. The notion of admissible graphs are defined in the same way as before. For each admissible graph, set
\[ c_\Gamma:= \int_{E_{k+2,n+1}} \bigwedge_{e\in E_\Gamma} d\theta_e.\]
Define an operator $H_0: C_*(A^\mathscr{W}) \ra \Omega_A^*[[t_1,\cdots,t_\mu]]$ by formula
\[ H_0(a_0|\cdots|a_n):= \sum_{\substack{\Gamma\in G'(k+2,n+1)\\ k\geq 0, l\geq 0\\ I=(i_1,\cdots,i_l)\in \left\{ 1,\cdots, N\right\}^j}} \frac{1}{k!} c_\Gamma \mathcal{U}_\Gamma(\partial_I,\phi_j,\underbrace{\mathscr{W},\cdots,\mathscr{W}}_{\mbox{$k$ copies}};a_0|\cdots|a_n) d\epsilon_I.\]
This operator is illustrated in the following diagram.
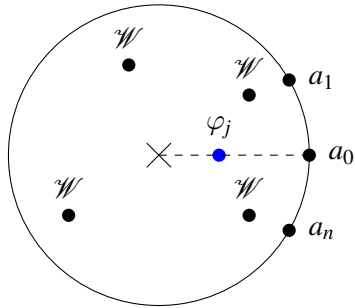
\[ \begin{tikzpicture}[baseline={(current bounding box.center)},scale=0.4]
\draw (0,0) circle (5);
\node[circle,fill=black,inner sep=0pt,minimum size=5pt,label=right:{$a_0$}] at (5,0) {};
\node[circle,fill=black,inner sep=0pt,minimum size=5pt,label=right:{$a_1$}] at (4.33,2.5) {};
\node[circle,fill=black,inner sep=0pt,minimum size=5pt,label=right:{$a_n$}] at (4.33,-2.5) {};
\draw (0,0) node[cross=5pt] {};
\node[circle,fill=blue,inner sep=0pt,minimum size=5pt,label=above:{$\varphi_j$}] at (2,0) {};
\node[circle,fill=black,inner sep=0pt,minimum size=5pt,label=above:{$\mathscr{W}$}] at (3,2) {};
\node[circle,fill=black,inner sep=0pt,minimum size=5pt,label=above:{$\mathscr{W}$}] at (-1,3) {};
\node[circle,fill=black,inner sep=0pt,minimum size=5pt,label=above:{$\mathscr{W}$}] at (-3,-2) {};
\node[circle,fill=black,inner sep=0pt,minimum size=5pt,label=above:{$\mathscr{W}$}] at (3,-2) {};
\draw[dashed] (0,0) to (5,0);
\end{tikzpicture}\]
Define another operator $H_1: C_*(A^\mathscr{W}) \ra \Omega_A^*[[t_1,\cdots,t_\mu]]$ in a similar way with the first marked point (corresponding to the point where $\phi_j$ is inserted) constraint between the framing of the central point $(0,0)$ and the boundary marked point $\overline{0}$. A typical configuration is illustrated in the following picture.
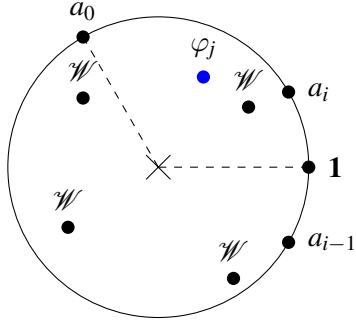
\[ \begin{tikzpicture}[baseline={(current bounding box.center)},scale=0.4]
\draw (0,0) circle (5);
\node[circle,fill=black,inner sep=0pt,minimum size=5pt,label=right:{$\bone$}] at (5,0) {};
\node[circle,fill=black,inner sep=0pt,minimum size=5pt,label=right:{$a_i$}] at (4.33,2.5) {};
\node[circle,fill=black,inner sep=0pt,minimum size=5pt,label=right:{$a_{i-1}$}] at (4.33,-2.5) {};
\node[circle,fill=black,inner sep=0pt,minimum size=5pt,label=above:{$a_0$}] at (-2.5,4.33) {};
\draw (0,0) node[cross=5pt] {};
\draw[dashed] (0,0) to (-2.5,4.33);
\draw[dashed] (0,0) to (5,0);
\node[circle,fill=blue,inner sep=0pt,minimum size=5pt,label=above:{$\varphi_j$}] at (1.5,3) {};
\node[circle,fill=black,inner sep=0pt,minimum size=5pt,label=above:{$\mathscr{W}$}] at (3,2) {};
\node[circle,fill=black,inner sep=0pt,minimum size=5pt,label=above:{$\mathscr{W}$}] at (-2.5,2.3) {};
\node[circle,fill=black,inner sep=0pt,minimum size=5pt,label=above:{$\mathscr{W}$}] at (-3,-2) {};
\node[circle,fill=black,inner sep=0pt,minimum size=5pt,label=above:{$\mathscr{W}$}] at (2.5,-3.7) {};
\end{tikzpicture}\]
\begin{Proposition}
For each $0\leq j\leq \mu$, the following identity holds:
\[ \big( \nabla_{\frac{\partial}{\partial t_j}}^{{\sf GM}}\pm \frac{1}{2} d_{DR}L_{\phi_j}\big)\circ \mathcal{U}_0^{Sh,\mathscr{W}} - \mathcal{U}_0^{Sh,\mathscr{W}}\circ \nabla^{{\sf Get}}_{\frac{\partial}{\partial t_j}} = u^{-1}\big( H\circ (b+uB) + (L_{\mathscr{W}}+ud_{DR})\circ H\big),\]
where $H=H_0+uH_1$.
\end{Proposition}

\begin{Proof}
Expanding both sides of the above equation, and moving the right hand side to the left, its $u^{-1}$-term is given by
\[ -\varphi_j\cap  \mathcal{U}_0^{Sh,\mathscr{W}} + \mathcal{U}_0^{Sh,\mathscr{W}} b^{1|1}(\frac{\partial \mu}{\partial t_j})-H_0 b-L_{\mathscr{W}} H_0.\]  
The vanishing of the above expression is proved by Calaque-Rossi~\cite{CalRos}. The $u^0$-term are given by
\[ \pm \frac{1}{2} d_{DR}L_{\phi_j} \mathcal{U}_0^{Sh,\mathscr{W}} +\mathcal{U}_0^{Sh,\mathscr{W}}B^{1|1}(\frac{\partial \mu}{\partial t_j})-H_1b-L_{\mathscr{W}}H_1+ [\frac{\partial}{\partial t_j}, \mathcal{U}_0^{Sh,\mathscr{W}}] -d_{DR}H_0-H_0B.\]
This vanishes as well, thanks to~\cite[Theorem 4.1]{CFW}. Note that here, we need to use the fact that
\[ [\frac{\partial}{\partial t_j}, \mathcal{U}_0^{Sh,\mathscr{W}}]= \mathcal{U}_1^{Sh,\mathscr{W}}(\phi_j).\] 
Finally, the $u^1$-term is given by
$H_1B+d_{DR}H_1$ which in fact vanishes separately, see~\cite[Paragraph 4.1]{CFW}.
\end{Proof}

Observe that since we assume that $W$ has an isolated singularity, cohomology classes of the twisted de Rham complex $H^*\big( \Omega^*_A\otimes \mathbb{K}[[t_1,\cdots,t_\mu]][[u]], L_{\mathscr{W}}+ud_{DR}\big)$ are represented by $d_{DR}$-closed elements. Thus the extra term $ \pm \frac{1}{2} d_{DR}L_{\phi_j} $ in the above Proposition will not appear, after taking cohomology. This yields the following

\begin{Proposition}~\label{prop:gm}
For each $0\leq j\leq \mu$, the following diagram is commutative.
\[\begin{CD}
HC^-_*(A^{\mathscr{W}}) @>\mathcal{U}_0^{Sh,\mathscr{W}} >> H^*\big( \Omega^*_A\otimes \mathbb{K}[[t_1,\cdots,t_\mu]][[u]], L_{\mathscr{W}}+ud_{DR}\big)\\
@V\nabla^{{\sf Get}}_{\frac{\partial}{\partial t_j}} VV   @VV \nabla_{\frac{\partial}{\partial t_j}}^{{\sf GM}} V\\
HC^-_*(A^{\mathscr{W}}) @>\mathcal{U}_0^{Sh,\mathscr{W}} >> H^*\big( \Omega^*_A\otimes \mathbb{K}[[t_1,\cdots,t_\mu]][[u]], L_{\mathscr{W}}+ud_{DR}\big)
\end{CD}\]
\end{Proposition}

\begin{Proposition}~\label{prop:uconn}
The following diagram is commutative.
\[\begin{CD}
HC^-_*(A^{\mathscr{W}}) @>\mathcal{U}_0^{Sh,\mathscr{W}} >> H^*\big( \Omega^*_A\otimes \mathbb{K}[[t_1,\cdots,t_\mu]][[u]], L_{\mathscr{W}}+ud_{DR}\big)\\
@V\nabla^{{\sf NC}}_{\frac{\partial}{\partial u}} VV   @VV \nabla_{\frac{\partial}{\partial u}} V\\
HC^-_*(A^{\mathscr{W}}) @>\mathcal{U}_0^{Sh,\mathscr{W}} >> H^*\big( \Omega^*_A\otimes \mathbb{K}[[t_1,\cdots,t_\mu]][[u]], L_{\mathscr{W}}+ud_{DR}\big)
\end{CD}\]
\end{Proposition}

\begin{proof}
Recall the formula of the connection operator $\nabla^{{\sf NC}}_{\frac{\partial}{\partial u}}$ from~\cite{CLT}:
\[ \nabla^{{\sf NC}}_{\frac{\partial}{\partial u}}:= \frac{\partial}{\partial u} + \frac{\Gamma}{2u}+ \frac{\iota(\prod_n (2-n) \mu_n )}{2u^2},\]
where $\Gamma(a_0|a_1|\cdots|a_n)=-n a_0|a_1|\cdots|a_n$, and $\iota(-)=b^{1|1}(-)+uB^{1|1}(-)$. The idea is to consider the following additional formal deformation in the $t$-direction:
\[ \mathscr{W}_t=\mathscr{W}+t\mathscr{W}.\]
Denote the corresponding family of $A_\infty$ strucutre by $\mu^t$. Its Kodaira-Spencer class in the $t$-direction is given by $\partial \mu^t/\partial t= \mathcal{U}_1^{{\mathscr{W}}}(\mathscr{W})=\sum_{k\geq 0} \frac{1}{k!}(1+t)^k \mathcal{U}_{k+1}(\mathscr{W}^{k+1})$. Its restriction to $t=0$ is given by 
\[\frac{\partial \mu^t}{\partial t}|_{t=0}=\sum_{k\geq 0} \frac{1}{k!} \mathcal{U}_{k+1}(\mathscr{W}^{k+1}).\] 
By definition, the $A_\infty$ structure on $A^{{\mathscr{W}}}$ is given by
\[ \mu_n = \sum_{k\geq 1} \frac{1}{k!} \mathcal{U}_k(\mathscr{W}^k)_n.\]
The subscript $n$ is the component of $\mathcal{U}_k(\mathscr{W}^k)$ that has $n$-inputs. We then compute the difference:
\begin{align*}
&   2\frac{\partial \mu^t}{\partial t}|_{t=0}-\prod_n (2-n)\mu_n\\
=& \sum_{k\geq 1, n \geq 1} (2k+n-2) \frac{1}{k!} \mathcal{U}_k(\mathscr{W}^k)_n\\
=& [\mu, T]_G \;\; \mbox{(The Gerstenhaber bracket)}
\end{align*}
The operator $T$ in the above formula is simply the exterior degree operator, i.e. $T(\epsilon_I)=|I|\epsilon_I$. The last equality follows from the fact that in order for the component $ \mathcal{U}_k(\mathscr{W}^k)_n$ to be non-vanishing, the number of edges of the admissible graphs must equal to the dimension of $C_{k,n}$ which is equal to $2k+n-2$. Since each edge differentiate once, the total operator reduces the exterior degree by exactly $2k+n-2$. 

Using the above formula, we compute
\begin{align*}
2u \nabla^{{\sf NC}}_{\frac{\partial}{\partial u}}=& 2u \frac{\partial}{\partial u} + \Gamma+ \frac{\iota(\prod_n (2-n) \mu_n )}{u}\\
=& \big( 2u \frac{\partial}{\partial u} + \Gamma +L_T +  \frac{\iota(2\frac{\partial \mu^t}{\partial t})}{u}\big)|_{t=0}\\
=&  \big( 2u \frac{\partial}{\partial u} + \Gamma +L_T + 2\frac{\partial}{\partial t}\big)|_{t=0} - \nabla^{{\sf Get}}_{2\frac{\partial}{\partial t}} |_{t=0}
\end{align*}
We have used Getzler's Cartan homotopy formula
\[ \iota(2\frac{\partial \mu^t}{\partial t}|_{t=0}) - \iota(-\prod_n (2-n)\mu_n) = u L_T,\]
in the second equality, with $L_T$ the Lie derivative action of Hochschild cochains on Hochschild chains. At this point, we apply the previous Corollary to obtain
\begin{align*}
& \mathcal{U}^{Sh,\mathscr{W}}_0 \circ 2u \nabla^{{\sf NC}}_{\frac{\partial}{\partial u}}\\
=& \mathcal{U}^{Sh,\mathscr{W}}_0 \big( 2u \frac{\partial}{\partial u} + \Gamma +L_T + 2\frac{\partial}{\partial t}\big)|_{t=0} -\big( \nabla^{{\sf GM}}_{2\frac{\partial}{\partial t}}\circ \mathcal{U}^{Sh,\mathscr{W}}_0 \big)|_{t=0}\\
=& \mathcal{U}^{Sh,\mathscr{W}}_0 \circ ( 2u \frac{\partial}{\partial u} ) +  \mathcal{U}^{Sh,\mathscr{W}}_0 \circ( \Gamma +L_T)- [2\frac{\partial}{\partial t},  \mathcal{U}^{Sh,\mathscr{W}}_0]|_{t=0}+\frac{2\iota_{\mathscr{W}}}{u}\circ  \mathcal{U}^{Sh,\mathscr{W}}_0
\end{align*}
Since $\mathcal{U}^{Sh,\mathscr{W}}_0$ is $u$-linear, it commutes with $2u\frac{\partial}{\partial u}$. For the remaining two terms, we observe that
\[ \mathcal{U}^{Sh,\mathscr{W}}_0 \circ( \Gamma +L_T)- [2\frac{\partial}{\partial t},  \mathcal{U}^{Sh,\mathscr{W}}_0]|_{t=0} =P\circ \mathcal{U}^{Sh,\mathscr{W}}_0,\]
where the operator $P: \Omega_A^*[[t_1,\cdots,t_\mu,u]] \rightarrow \Omega_A^*[[t_1,\cdots,t_\mu,u]] $ is defined by $$ P(\epsilon_I (d\epsilon_1)^{k_1}\cdots(d\epsilon_N)^{k_N} )= |I|\cdot \epsilon_I (d\epsilon_1)^{k_1}\cdots(d\epsilon_N)^{k_N}$$ extended $\mathbb{K}[[t_1,\cdots,t_\mu,u]]$-linearly. That is, it simply records the tensor degree of the coefficient of a differential form. The above identity follows from the fact that the dimension of the configuration space $D_{k+1,n+1}$ is $2k+n$, and hence one needs to contract $2k+n$ tensor degrees for $\mathcal{U}_k^{Sh}(\mathscr{W},\cdots,\mathscr{W};a_0|a_1|\cdots|a_n)$ to be non-vanishing. Putting the above identities together, we obtain
\[ \mathcal{U}^{Sh,\mathscr{W}}_0 \circ 2u \nabla^{{\sf NC}}_{\frac{\partial}{\partial u}}=(2u\frac{\partial}{\partial u} + P +\frac{2\iota_{\mathscr{W}}}{u})\circ \mathcal{U}^{Sh,\mathscr{W}}_0.\]
Finally, observe that the $P$ operator corresponds to the differential form degree operator on $\Omega_S^*$, and since we are assuming that $W$ has an isolated singularity at origin, the cohomology $H^*\big( \Omega_S^*[[t_1,\cdots,t_\mu,u]], d\mathscr{W}+ud_{DR}\big)$ is supported at the top differential form degree $N$. Hence the operator $P\equiv N\cdot \id$ after taking cohomology. Thus, the proposition is proved.
\end{proof}

\end{document}